\documentclass[12pt,reqno]{amsart}

\makeatletter
\def\@settitle{%
	\begin{center}%
		\normalfont\Large\bfseries 
		\@title
	\end{center}%
}
\makeatother

\usepackage{mathrsfs}
\usepackage{enumitem}   
\usepackage{graphicx}
\usepackage[colorlinks = true,
linkcolor = blue,
urlcolor  = blue,
citecolor = blue,
anchorcolor = blue]{hyperref}

\newtheorem{theorem}{Theorem}

\newtheorem{proposition}
{Proposition}
\newtheorem{corollary}
{Corollary}
\newtheorem{remark}{Remark}

\newtheorem{lemma}{Lemma}

\newfont{\bb}{msbm10 at 12pt}

\def\SB{\mathbf{S}_g}
\def\TSB{\mathcal{S}_g}

\def\TSBH{\mathcal{S}_b}

\def\EH{\mathcal{EK}_b}

\def\<{\langle}     
\def\>{\rangle}     
\def\div{{\rm div}}

\def\Tr{{\rm tr}}

\newcommand{\bal}{\begin{align}}      \newcommand{\eal}{\end{align}}
\newcommand{\ba}{\begin{array}}      \newcommand{\ea}{\end{array}}
\newcommand{\bc}{\begin{center}}     \newcommand{\ec}{\end{center}}
\newcommand{\be}{\begin{enumerate}}  \newcommand{\ee}{\end{enumerate}}
\newcommand{\beq}{\begin{eqnarray}}  \newcommand{\eeq}{\end{eqnarray}}
\newcommand{\beQ}{\begin{eqnarray*}} \newcommand{\eeQ}{\end{eqnarray*}}
\newcommand{\bi}{\begin{itemize}}    \newcommand{\ei}{\end{itemize}}
\newcommand{\bt}{\begin{tabular}}    \newcommand{\et}{\end{tabular}}
\newcommand{\bdm}{\begin{displaymath}} \newcommand{\edm}{\end{displaymath}}

\newcommand{\D}{D\!\!\!\!/\,}

\newcommand{\ETSB}{\mathcal{S}\!\!\!\!/\,}
\newcommand{\nb}{\nabla\!\!\!\!/\,}

\def\qed{\hfill{q.e.d.}\smallskip\smallskip}

\begin{document}
	
	\title[]{Positive energy-momentum theorems for asymptotically  AdS spin initial data sets with charge}       
	
	\author{Simon Raulot}
	\address[Simon Raulot]{Univ Rouen Normandie, CNRS, Normandie Univ, LMRS UMR 6085, F-76000 Rouen, France}
	\email{simon.raulot@univ-rouen.fr}
	
	\begin{abstract}
		For complete spin initial data sets with an asymptotically anti--de Sitter end, we introduce a charged energy--momentum defined as a linear functional arising from the Einstein--Maxwell constraints. Under a dominant energy condition adapted to the presence of a negative cosmological constant, we establish positive energy--momentum theorems, showing in particular that this functional is non--negative on a natural real cone. We place particular emphasis on the case where the manifold carries a compact inner boundary. In the time--symmetric setting, this yields a mass--charge inequality for asymptotically hyperbolic manifolds with charge.
	\end{abstract}

	
	
	\thanks{}
	
	\date{\today}   
	
	\maketitle 
	\pagenumbering{arabic}	
	
	
	\section{Introduction}
	
	
	The positive mass theorem is one of the cornerstones of mathematical general relativity: under suitable energy conditions, the total energy of an isolated gravitating system is non--negative, with rigidity characterizing the relevant ground state. In the asymptotically flat setting, this paradigm goes back to the seminal works of Schoen--Yau \cite{SchoenYau1,SchoenYau6} and Witten \cite{Witten1}, and has since inspired a vast literature on mass--type invariants, boundary effects, and rigidity phenomena.
	
	In the presence of a negative cosmological constant, the natural vacuum background becomes anti--de Sitter spacetime. Motivated both by geometric questions (e.g.\ scalar curvature rigidity) and by developments in high--energy physics (such as the AdS/CFT correspondence and related notions of conserved charges), a theory of total energy--momentum for asymptotically AdS systems has been developed over the last decades. In the time--symmetric case, positivity and rigidity for asymptotically hyperbolic manifolds were obtained by Wang \cite{Wangx3} and by Chru\'sciel--Herzlich \cite{ChruscielHerzlich}. 
	
	Beyond the time--symmetric regime, one is led to consider general initial data sets $(M^n,g,K)$. For asymptotically AdS initial data, the energy--momentum and additional conserved quantities associated with the AdS isometry group were investigated by Chru\'sciel--Nagy \cite{ChruscielNagy}, and further developed in subsequent works, including spinorial proofs of positivity in the AdS--asymptotically hyperbolic setting by Maerten \cite{Maerten2}, as well as rigidity and sharp inequalities involving angular momentum and center of mass by Chru\'sciel--Maerten--Tod \cite{ChruscielMaertenTod}. 
	
	The purpose of the present paper is to incorporate an \emph{electric field} into this picture and to establish \emph{positive energy--momentum theorems with charge} for asymptotically AdS initial data sets in arbitrary dimension $n\ge 3$ under a dominant energy condition adapted to the Einstein--Maxwell constraints. While electric charge and the associated conserved quantities are well known in the physics literature, our aim here is to place them within a rigorous geometric framework in the setting of asymptotically AdS initial data. Concretely, we consider complete charged initial data sets $(M^n,g,K,E)$ possessing at least one asymptotically AdS end, modeled on the hyperbolic metric, and we define a \emph{charged energy--momentum} as a linear functional on the kernel of the adjoint of a suitable linearized constraint operator, following Michel's formalism \cite{Michel}. This charged energy--momentum extends the (uncharged) AdS energy--momentum of Chru\'sciel--Nagy \cite{ChruscielNagy} and contains, as a special component, the total electric charge.
	
	Our main results rely on a spinorial approach which can be viewed as a natural continuation of earlier works in dimension three by Bartnik--Chru\'sciel \cite{BartnikChrusciel}, as well as of the author's previous results in the asymptotically flat setting \cite{Raulot16}. The key analytic input is a Witten--type
	argument based on a modified Dirac--Witten operator and on the existence, on the hyperbolic background, of \emph{extrinsic imaginary Killing spinors}. These spinors
	are naturally associated with Lorentzian imaginary Killing spinors on AdS, restricted to the totally geodesic hyperbolic slice, and they generate a distinguished positivity cone in the space of Killing initial data at infinity. The charged energy--momentum is shown to be non--negative on this cone under the dominant energy condition, yielding a positivity statement for the charged conserved quantities (see Theorem \ref{PET-Spherical}). We also consider initial data sets with a compact inner boundary, leading to positive energy--momentum theorems (see Theorems \ref{PositiveETS-Boundary} and \ref{PositiveETS-Boundary-APS}) and, in particular, to a new such theorem in the uncharged asymptotically anti--de Sitter setting (see Corollary \ref{PositiveSTB}). Our main results may be viewed as charged analogues of the positivity results obtained by Maerten in the asymptotically hyperbolic setting \cite{Maerten1,Maerten2}.
	
	In the time--symmetric regime, the AdS energy--momentum reduces to the hyperbolic mass functional, and our method yields a mass--charge inequality for complete charged data with an asymptotically hyperbolic end (see Theorem \ref{PMTC-Hyp}), together with a corresponding statement in the presence of an inner boundary (see Theorem~\ref{PMTCB-Hyp-TS}, where two alternative boundary conditions are studied). 
	
	Finally, let us stress that the restriction to the asymptotically AdS setting with spherical conformal infinity is deliberate. While more general notions of asymptotically locally AdS ends are available and play an important role in other contexts, the spinorial method developed here does not directly extend to such
	settings. In particular, for non--spherical conformal infinities, such as toroidal ones, the extrinsic imaginary Killing spinors available at infinity do not capture the contribution of the electric charge in the associated energy--momentum expressions. As a consequence, the present approach does not yield a mass--charge inequality in those cases. For broader perspectives on AdS energy definitions and positivity mechanisms, including links with holography and BPS--type inequalities, we refer for instance to \cite{Rallabhandi} and the references therein.
	
	We do not address the rigidity or equality cases in the present work; a detailed analysis of these questions will be carried out in a forthcoming paper \cite{Raulot17}.
	
	\section{Definitions and statement of the results}
	
	\subsection{Asymptotically AdS initial data sets with charge}\label{AdSIDS}
	
	An \emph{initial data set with charge} (or a \emph{charged initial data set}) $(M^n,g,K,E)$ is a Riemannian $n$-dimensional manifold $(M^n,g)$ e\-quip\-ped with a symmetric $(0,2)$-tensor $K$ and a vector field $E$ on $M$. In this work, such a quadruplet is said to have an \emph{asymptotically AdS end} if the manifold contains a region $M_{ext}\subset M$ together with a diffeomorphism $\Psi^{-1}:M_{ext}\rightarrow\mathcal{E}_R:=(R,\infty)\times\mathbb{S}^{n-1}$ with $R>0$, called a chart at infinity, such that, in this chart, we have
	\begin{equation}\label{AsymptoticDecay}
		e_g:= g-b=\mathcal{O}_1(e^{-\tau r}), \quad  K=\mathcal{O}(e^{-\tau r}),\quad  E=\mathcal{O}(e^{-\tau r})
	\end{equation}
	for some $\tau>n/2$. Here a tensor $T$ on $M_{ext}$ satisfies $T=\mathcal{O}_ k(e^{-\tau r})$ for $k\in\mathbb{N}$ if
	\begin{align*}
		|T|_b+|\nabla^b T|_b+\cdots+|\underbrace{\nabla^b\cdots\nabla^b}_{k\text{ factors}} T|_b\leq Ce^{-\tau r}
	\end{align*}
	as $r$ goes to infinity, for some constant $C>0$, and where $|\,\cdot\,|_b$ and $\nabla^b$ denote respectively the norm of a tensor field and the covariant derivative with respect to the background metric $b$. The background metric appearing in \eqref{AsymptoticDecay} is the hyperbolic metric with constant scalar curvature $-n(n-1)$ given by
	\begin{align*}
		b = dr^2 + \sinh^2 (r)\,\check h,
	\end{align*}
	where $\check h$ is the standard metric on the $(n-1)$-dimensional sphere $\mathbb{S}^{n-1}$ of constant sectional curvature equal to $1$. In this situation, $(\mathbb{H}^n,b)$ provides a time-symmetric initial data set with an asymptotically AdS end (with $E=0$) and gives rise to a static solution of the Einstein equations with negative cosmological constant. The notion of asymptotically AdS end adopted in this work is based on the standard hyperbolic model. A broader notion of asymptotically locally AdS ends, allowing different geometries at infinity, is available in the literature (see, e.g., \cite{ChruscielNagy}). In the time-symmetric setting, that is when $K=0$, the triplet $(M^n,g,E)$ is said to have an \emph{asymptotically hyperbolic end}. 
	
	\medskip
	
	A natural and physically relevant class of initial data sets with charge arises from electrovacuum spacetimes with negative cosmological constant. Let $(\mathcal{M}^{n+1},\mathbf{g})$ be a time-oriented Lorentzian manifold and let $F$ be a two-form, usually called the Faraday two-form. The Einstein--Maxwell equations with cosmological constant $\Lambda<0$ are
	\begin{align*}
		\mathrm{Ric}_{\mathbf g} - \tfrac12 R_{\mathbf g}\,\mathbf g + \Lambda\,\mathbf g = \mathcal{T},
		\qquad dF = 0, 
		\qquad d(\star_{\mathbf g}F)=0
	\end{align*}
	where $\mathrm{Ric}_{\mathbf g}$ and $R_{\mathbf g}$ denote respectively the Ricci and scalar curvatures of ${\mathbf g}$.  The stress--energy tensor of the electromagnetic field is
	\begin{align*}
		\mathcal{T}_{\alpha\beta}
		=2\Big(F_{\alpha\mu}F_{\beta}{}^{\mu}
		-\tfrac14 \mathbf{g}_{\alpha\beta}\,F_{\mu\nu}F^{\mu\nu}\Big).
	\end{align*}
	We normalize $\Lambda=-n(n-1)/2$. Let $M^n\subset\mathcal{M}^{n+1}$ be a spacelike hypersurface with induced metric $g$, unit future-directed normal $N$, and second fundamental form $K$. The Gauss and Codazzi equations imply that the induced initial data satisfy
	\begin{align*}
		R_g + (\Tr_g K)^2 - |K|_g^2 = 2\mathcal{T}(N,N)-n(n-1),
		\qquad
		\mathrm{div}_g K - d(\Tr_g K) = -\,\mathcal{T}(N,\cdot)^\sharp
	\end{align*}
	with $R_g$ the scalar curvature of $g$. In the purely electric Einstein--Maxwell case,
	\begin{align*}
		F = c_n\, N^\flat\wedge E^\flat,
		\qquad c_n=\sqrt{\tfrac{(n-1)(n-2)}{2}},
	\end{align*}
	where $E$ is tangent to $M$. A straightforward computation shows that
	\begin{align*}
		\mathcal{T}(N,N)=\tfrac12 (n-1)(n-2)|E|_g^2, \qquad
		\mathcal{T}(N,\cdot)=0.
	\end{align*}
	Therefore, the quadruplet $(M^n,g,K,E)$ defines an initial data set with charge. Note that the Maxwell equation $d(\star_{\mathbf g}F)=0$ further reduces to $\mathrm{div}_g E=0$. Many explicit electrovacuum solutions with $\Lambda<0$, such as the Reissner--Nordström--AdS family, give rise to complete charged initial data sets with an asymptotically AdS end in the sense defined above. 
	
	\subsection{The charged energy-momentum}\label{CEM}
	
	In order to provide a geometric definition of the energy–momentum in this context, we follow the framework developed by Michel \cite{Michel}. We are therefore led to consider the following constraints map
	$$
	\begin{array}{rcl}
		\Phi: \mathcal{M} \times\Gamma\left(S^2M\right)\times\Gamma\left(TM\right) & \longrightarrow & C^\infty(M)\times\Gamma(T^\ast M)\times C^\infty(M) \\[4pt]
		\phantom{f:}\, h:=(g,K,E) & \longmapsto & \big(\Phi_1(g,K,E),\Phi_2(g,K,E),\Phi_3(g,K,E)\big)
	\end{array}
	$$	
	where 
	$$
	\left\lbrace
	\begin{array}{lll}
		\Phi_1(g,K,E) & = & R_g+(\Tr_g K)^2-|K|^2_g -(n-1)(n-2)|E|^2_g,\\[2pt]
		\Phi_2(g,K,E) & = & 2\left(\div_g K-d(\Tr_g K)\right),\\[2pt]
		\Phi_3(g,K,E)& = & 2(n-1)\div_g E,
	\end{array}
	\right.
	$$
	and where $\mathcal{M}$ denotes the set of Riemannian metrics on the manifold $M$ and $S^2M$ is the space of covariant symmetric $2$-tensors. Let $h_0:=(b,0,0)$ and $e:= h-h_0$. One computes that 
	\begin{align*}
		\< \eta,\Phi( h)-\Phi(h_0)\>_b =\div_b \mathbb{U}(\eta,e)+Q(\eta,e)
	\end{align*} 
	for $\eta\in{\rm Ker} \left(D_{h_0}\Phi\right)^\ast$ and where $Q(\eta,e)$ is the remainder term which is linear in $\eta$ and at least quadratic in $e$. Here the integrand $\mathbb{U}(\eta,e)$ is the $1$-form appearing in the following integration by parts formula:
	\begin{align*}
		\<\eta,D_{h_0}\Phi(e)\>_b=\div_b\mathbb{U}(\eta,e)+\<(D_{h_0}\Phi)^\ast\eta,e\>_b 
	\end{align*}
	whose last term will be rewritten as 
	\begin{align*}
		\<(D_{h_0}\Phi)^\ast\eta,e\>_b =\<(D_{h_0}\Phi)^\ast_g(\eta),e_g\>_b+\<(D_{h_0}\Phi)^\ast_K(\eta), K\>_b+\<(D_{h_0}\Phi)^\ast_E(\eta), E\>_b
	\end{align*}
	for all $\eta=(V,\alpha,f)\in C^\infty(M)\times\Gamma(T^\ast M)\times C^\infty(M)$. 
	However, the  linearized operators at $h_0$ are
	$$\left\lbrace
	\begin{array}{rcl}
		D_{h_0}\Phi_1(e) & = & \div_{b}\left(\div_{b}e_g-d(\Tr_be_g)\right)+(n-1)\Tr_be_g\\
		D_{h_0}\Phi_2(e) & = & 2\big(\div_b  K-d(\Tr_b  K)\big)\\
		D_{h_0}\Phi_3(e) & = & 2(n-1)\div_b E
	\end{array}
	\right.
	$$
	and so one computes that
	\begin{equation}\label{AdjLinOp}
		\left\lbrace
		\begin{array}{lll}
			(D_{h_0}\Phi)^\ast_g(V,\alpha,f) & = & \nabla^bdV-(\Delta_bV)b+(n-1)Vb \\
			(D_{h_0}\Phi)^\ast_K(V,\alpha,f) & = & -\mathcal{L}_{\alpha^\sharp} b+2(\div_b\alpha)b\\
			(D_{h_0}\Phi)^\ast_E(V,\alpha,f) & = & -2(n-1)\nabla^bf \\
		\end{array}
		\right.
	\end{equation}
	and $\mathbb{U}(\eta,e)=\mathbb{U}_1(V)+\mathbb{U}_2(\alpha)+\mathbb{U}_3(f)$ where 
	\begin{equation}
		\left\lbrace
		\begin{array}{lll}
			\mathbb{U}_1(V) & = & V\left(\div_b e_g -d(\Tr_b e_g)\right)-i_{\nabla^b V}e_g+(\Tr_b e_g)dV,\nonumber\\[2pt]
			\mathbb{U}_2(\alpha) & = & 2\left(i_{\alpha^{\sharp}} K-(\Tr_b K)\alpha\right),\\[2pt]
			\mathbb{U}_3(f) & = & 2(n-1)f E^\flat .
		\end{array}
		\right.
	\end{equation}
	One then observes that if the decay assumptions (\ref{AsymptoticDecay}) hold for $h=(g,K,E)$ and
	\begin{align}\label{IntegrabilityCondition} 
		\<\eta,\Phi(h)-\Phi(h_0)\>_b\in L^1\big(\mathcal{E}_R,d\mu_b\big)
	\end{align}
	for all $\eta\in {\rm Ker} \left(D_{h_0}\Phi\right)^\ast$ then the total charge 
	\begin{eqnarray*}
		m(h,\Psi,\eta):=\lim_{r\rightarrow\infty}\oint_{S_r}\mathbb{U}(\eta,e)(\nu_{b_r})d\sigma_{b_r}
	\end{eqnarray*}
	is finite. Here the right-hand side is computed in the asymptotically AdS coordinates given by $\Psi$ on the distinguished end. Moreover $S_r$ denotes the slices $\{r\}\times\mathbb{S}^{n-1}$ in the product $(R,\infty)\times\mathbb{S}^{n-1}$, $b_r$ its metric induced by $b$ and $\nu_{b_r}$ its unit normal pointing toward infinity with respect to $b$. Furthermore, as shown in \cite{Michel}, although this definition seems to depend on the chart $\Psi$, it is geometrically well defined whenever (\ref{AsymptoticDecay}) and (\ref{IntegrabilityCondition}) hold. 
	
	From the expressions (\ref{AdjLinOp}), it follows that the kernel of the adjoint of the linearized operator of $\Phi$ is isomorphic to the vector space 
	\begin{align*}
		\mathfrak{Kill}(AdS^{n,1}) \oplus \mathbb{R} \simeq \mathfrak{so}(n,2) \oplus\mathbb{R} \simeq \mathcal{N}_b \oplus \mathfrak{Kill}(\mathbb{H}^n) \oplus \mathbb{R},
	\end{align*}
	where $\mathcal{N}_b:=\left\{V\in C^\infty(\mathbb{H}^n) \mid \nabla^{b}dV=Vb\right\}$ is the vector space of hyperbolic static potentials, $\mathfrak{Kill}(AdS^{n,1})$ and $\mathfrak{Kill}(\mathbb{H}^n)\simeq\mathfrak{so}(n,1)$ denote the Lie algebra of Killing vector fields on $AdS^{n,1}$ and on $\mathbb{H}^n$. Subsequently,  the \emph{charged energy–mo\-men\-tum} $\Xi$ of a charged initial data set $(M^n,g,K,E)$ is defined as the linear form on the kernel of the adjoint of $D_{h_0}\Phi$ by 
	\begin{align*}
		\Xi(V,\alpha,f):=\frac{1}{2(n-1)\omega_{n-1}}\lim_{r\rightarrow\infty}\oint_{S_r}\Big(\mathbb{U}_1(V)+\mathbb{U}_2(\alpha)+\mathbb{U}_3(f)\Big)(\nu_{b_r})d\sigma_{b_r}
	\end{align*}
	where $\omega_{n-1}$ denotes the volume of the unit sphere $(\mathbb{S}^{n-1},\check h)$. It contains, as special cases, the classical (uncharged) energy--momentum defined by Chru\'sciel and Nagy \cite{ChruscielNagy}, obtained by restricting
	$\Xi$ to
	\begin{align*}
		(V,\alpha)\in \mathcal{N}_b\oplus\mathfrak{Kill}(\mathbb{H}^n)
		\longmapsto \Xi(V,\alpha,0),
	\end{align*}
	and the \emph{electric charge}, given by
	\begin{align*}
		Q:=\Xi(0,0,1)
		= \frac{1}{\omega_{n-1}}\lim_{r\to\infty}
		\oint_{S_r} E^\flat(\nu_{b_r})\, d\sigma_{b_r}.
	\end{align*}
	\begin{remark}
		Finally, it is worth noting that this definition extends without difficulty to the broader class of asymptotically locally AdS ends considered in \cite{ChruscielNagy} (see \cite{ChruscielHerzlich} for the time-symmetric setting).
	\end{remark}
	
	\subsection{The dominant energy condition}\label{DEC}
	
	In order to ensure some positivity of the charged energy-momentum, we need a local positivity hypothesis, usually referred to as the dominant energy condition. More precisely, the initial data set with charge satisfies the \emph{dominant energy condition} if the vector $\Phi( h)-\Phi(h_0)\in\mathbb{R}^{n+1,1}$ is a future-directed causal vector for a given triplet $h=(g,K,E)$ on $M$. Equivalently, if 
	\begin{align*}
		\mu:=\frac{1}{2}\Phi_1(h),\quad J:=\frac{1}{2}\Phi_2(h),\quad \varpi:=\frac{1}{2}\Phi_3(h)
	\end{align*}
	denote respectively the associated energy, current and charge densities, the dominant energy condition then reads
	\begin{align}\label{DEC-General}
		\mu+\frac{n(n-1)}{2}\geq\sqrt{|J|_g^2+\varpi^2}.
	\end{align}
	Spacelike hypersurfaces arising from solutions of the Einstein-Maxwell equations satisfy this condition since $\Phi(h)-\Phi(h_0)=0$ in that case. 
	
	\subsection{Positive charged energy-momentum theorems}\label{CEMT}
	
	Before we can state our main results, it is necessary to fix some notation and recall a few preliminary constructions related to the background geometry and the spinorial formalism (see Section \ref{IntroSpinors} for the precise definitions).
	
	\medskip
	
	Our assumptions on the background metric ensure the existence of an \emph{extrinsic imaginary Killing spinor} for the metric $b$, namely a section $\varphi\in\Gamma(\TSBH)$ of a certain Dirac bundle over $\mathbb{H}^n$ such that
	\begin{align}\label{EKS}
		\forall X \in \Gamma(T\mathbb{H}^n), \quad\nabla^b_X \varphi = \mp \frac{i}{2} c_b(X)\varphi.
	\end{align}
	The terminology \emph{extrinsic imaginary Killing spinors} reflects the fact that such spinor fields are obtained by restriction of Lorentzian imaginary Killing spinors on anti--de Sitter spacetime to the hyperbolic space, viewed as a totally geodesic spacelike hypersurface. In the following, we denote by $\EH^\pm$ the space of such spinors. To each such spinor field, one associates via the map 
	\begin{align*}
		\zeta \in \EH^\pm \mapsto (V_\zeta,\alpha_\zeta,f_\zeta)
	\end{align*}
	an element of the kernel of $\left(D_{h_0}\Phi\right)^\ast$ where
	\begin{equation}
		\left\lbrace
		\begin{array}{ccl}
			V_\zeta & = & |\zeta|^2_b,\\
			\alpha_{\zeta}(X) & = & \<c_b(X)\gamma_b\zeta, \zeta \>_b, \quad\forall X\in\Gamma(T\mathbb{H}^n),\nonumber\\
			f_{\zeta} & = & \<\gamma_b\zeta,\zeta\>_b. 
		\end{array}
		\right.
	\end{equation}
	The image of this quadratic map is a real cone, denoted by $\mathscr{C}_b$; it can be parametrized by a map $\mathcal{K}:\Sigma_n\rightarrow \mathscr{C}_b$ defined on a finite--dimensional complex vector space (see Section \ref{PMR} for more details). We can now state the charged analogue of Maerten's theorem \cite[Theorem 1.3]{Maerten2}.
	\begin{theorem}\label{PET-Spherical}
		Let $(M^n,g,K,E)$, $n\geq 3$, be a complete spin initial data set with charge containing at least one asymptotically AdS end. Assume that the integrability condition (\ref{IntegrabilityCondition}) holds, as well as the dominant energy condition (\ref{DEC-General}). Then the linear form $\Xi\circ\mathcal{K}$ is non-negative. 
	\end{theorem}
	
	\medskip
	
	We also consider charged initial data sets with an asymptotically AdS end and carrying a compact inner boundary. In this situation, the strength of the gravitational field in a neighborhood of an $(n-1)$-hypersurface  $\Sigma$ can be characterized by its null expansions 
	\begin{align}\label{FuturePastTrapped}
		\theta_\pm := H \pm \Tr_\Sigma K 
	\end{align}
	where $H$ denotes the mean curvature with respect to the unit normal $\nu$ (pointing towards spatial infinity) and $\Tr_\Sigma K$ is the trace along $\Sigma$ of the tensor $K$. The null expansions measure the rate of change of the  area of a congruence of light rays emitted by the surface, in the outward future direction (with $\theta_+$) and outward past direction (with $\theta_-$). Thus the gravitational field is interpreted as being strong near $\Sigma$ if $\theta_+\leq 0$ (resp. $\theta_-\leq 0$); in that
	case $\Sigma$ is referred to as a future (resp. past) trapped surface. Future or past 
	apparent horizons arise as boundaries of future or past trapped regions and satisfy the equation $\theta_+=0$ (resp. $\theta_-=0$). We can then state the following result.
	\begin{theorem}\label{PositiveETS-Boundary}
		Under the assumptions of Theorem \ref{PET-Spherical}, suppose that $(M^n,g,K,E)$ has a compact inner boundary $\partial M$, each of whose connected components is either a future or a past trapped surface. Then the conclusion of Theorem \ref{PET-Spherical} remains valid.
	\end{theorem}
	
	Beyond trapped surface conditions, positivity results for initial data sets with a compact inner boundary can also be formulated under geometric constraints involving both the extrinsic geometry of the boundary and its conformal class. In this setting, the boundary contribution to the spinorial identities naturally leads to the appearance of the scalar quantity
	\begin{align*}
		\mathcal{H}_{max} := H+\sqrt{(\Tr_{\partial M} K)^2+|K(\nu)^\top|^2_{\mathfrak{g}}+(n-1)^2E_\nu^2},
	\end{align*}
	where $\mathfrak{g}:=g_{|\partial M}$, $K(\nu)^\top$ denotes the tangential part of $K(\nu)$, and $E_\nu$ is the normal component of the electric field. This combination arises from the zeroth order boundary term associated with the modified Dirac operator introduced in Section~\ref{BoundaryPart1}. The size of $\mathcal{H}_{max}$ is controlled in terms of the conformal geometry of the boundary $(\partial M,[\mathfrak{g}])$ through its Yamabe constant, defined by
	\begin{align*}
		\mathcal{Y}(\partial M,[\mathfrak{g}]):= \inf_{\overline{\mathfrak{g}}\in[\mathfrak{g}]} \left\{
		\frac{\oint_{\partial M}R_{\overline{\mathfrak{g}}}\,
			d\sigma_{\overline{\mathfrak{g}}}}
		{{\rm Vol}(\partial M,\overline{\mathfrak{g}})^{\frac{n-3}{n-1}}}
		\right\}
	\end{align*}
	where $R_{\overline{\mathfrak{g}}}$ and ${\rm Vol}(\partial M,\overline{\mathfrak{g}})$ denote respectively the scalar curvature and the volume of $\partial M$ computed with respect to a metric $\overline{\mathfrak{g}}$ in the conformal class $[\mathfrak{g}]$. Then we prove the following version of the positive-energy momentum theorem: 
	\begin{theorem}\label{PositiveETS-Boundary-APS}
		Under the assumptions of Theorem~\ref{PET-Spherical}, suppose that $(M^n,g,K,E)$ has a compact inner boundary $\partial M$ with positive Yamabe
		invariant satisfying
		\begin{align}\label{YamabeHmax}
			\frac{1}{n-1}\mathcal{H}_{max}\leq\sqrt{\frac{1}{(n-1)(n-2)} \mathcal{Y}(\partial M,[\mathfrak{g}]){\rm Vol}(\partial M,\mathfrak{g})^{-\frac{2}{n-1}}+1}.
		\end{align}
		Then the conclusion of Theorem~\ref{PET-Spherical} remains valid.
	\end{theorem}
	If the boundary is not connected, we assume that each connected component of $\partial M$ has positive Yamabe invariant and satisfies \eqref{YamabeHmax}. In dimension $n=3$, the Gauss--Bonnet theorem implies that $\partial M$ is a $2$--sphere, and the condition \eqref{YamabeHmax} reduces to
	\begin{align*}
		\frac{1}{2}\mathcal{H}_{max}\leq
		\sqrt{\frac{4\pi}{{\rm Vol}(\partial M,\mathfrak{g})}+1}.
	\end{align*}
	In this form, the above theorem yields a sharp positive energy theorem for (uncharged) asymptotically AdS initial data sets with a compact inner boundary, whose equality case is realized by the exterior of a geodesic ball in hyperbolic space.
	\begin{corollary}\label{PositiveSTB}
		Under the assumptions of Theorem~\ref{PET-Spherical} with $E=0$, suppose that $(M^n,g,K)$ has a compact inner boundary $\partial M$ with positive Yamabe
		invariant satisfying
		\begin{align*}
			\frac{1}{2}\Big(H+\sqrt{(\Tr_{\partial M} K)^2+|K(\nu)^\top|^2_{\mathfrak{g}}}\Big)\leq\sqrt{\frac{1}{(n-1)(n-2)} \mathcal{Y}(\partial M,[\mathfrak{g}]){\rm Vol}(\partial M,\mathfrak{g})^{-\frac{2}{n-1}}+1}.
		\end{align*}
		Then the conclusion of Theorem~\ref{PET-Spherical} remains valid.
	\end{corollary}
	
	\medskip
	
	In the time-symmetric setting, the uncharged energy-momentum reduces to a linear form on $\mathcal{N}_b$. As indicated in the previous section, the space $\mathcal{N}_b$ corresponds to the $(n+1)$-dimensional vector space of hyperbolic static potentials, with the following basis of functions:
	\begin{align*}
		V_{(0)}=\cosh r,\;\;V_{(1)}=x^1\sinh r,\;\;..., \;\;V_{(n)}=x^n\sinh r,
	\end{align*}
	where $x^1,\cdots,x^n$ are the coordinate functions on $\mathbb{R}^n$ restricted to $\mathbb{S}^{n-1}$. In the upper sheet of the unit hyperboloid model of $\mathbb{H}^n$, the functions $V_{(i)}$ are the restrictions to $\mathbb{H}^n$ of the coordinate functions on the Minkowski spacetime $\mathbb{R}^{n,1}$. The space $\mathcal{N}_b$ is naturally endowed with a Lorentzian metric $\eta$ of signature $(+,-,\cdots,-)$, induced by the action of the group of isometries $O^+(1, n)$ of the hyperbolic metric. This metric is characterized by the condition that the above basis is orthonormal, that is
	\begin{align*}
		\eta(V_{(0)},V_{(0)})=1\quad\text{and}\quad\eta(V_{(i)},V_{(i)})=-1
	\end{align*}
	for all $i=1,\cdots,n$. A time-orientation on $\mathcal{N}_b$ is provided by declaring that a time-like vector $X=X^{(\mu)}V_{(\mu)}\in\mathcal{N}_b$ is future-directed if $X^{(0)}>0$. It follows that if we let
	\begin{align}\label{MassCoeffSpherical}
		m_{(\mu)}:=\Xi(V_{(\mu)},0,0),\quad\mu=0,\cdots,n, 
	\end{align}
	then the real number $m$ defined (up to a sign) by
	\begin{align*}
		m^2:=\eta({\bf m},{\bf m})=\Big|m_{(0)}^2-\sum_{j=1}^nm_{(j)}^2\Big|
	\end{align*}
	is called the {\it mass} of $(M^n,g)$. Here ${\bf m}$ is the co-vector with components $m_{(\mu)}$. It is important to point out that the causal character as well as the future/past pointing nature of ${\bf m}$ are geometric invariants. Then if ${\bf m}$ is causal and future-directed, it is natural to choose $m:=\sqrt{\eta({\bf m},{\bf m})}$. This property is ensured under \eqref{IntegrabilityCondition} and \eqref{DEC-General}, which in the time-symmetric case reduce respectively to 
	\begin{align}\label{IntegrabilityHyperbolic}
		e^r\big(R_g+n(n-1)\big)\in L^1\big(\mathcal{E}_R,d\mu_b\big)
	\end{align}
	and $R_g\geq -n(n-1)$, by the classical hyperbolic positive mass theorem obtained by \cite{Wangx3} and Chru\'sciel and Herzlich \cite{ChruscielHerzlich}. In the time-symmetric charged setting, the dominant energy condition
	\eqref{DEC-General} reduces to 
	\begin{align}\label{DEC-TS}
		R_g\geq -n(n-1)+(n-1)(n-2)|E|^2_g+ 2(n-1)|\div_g E|
	\end{align}
	which already implies that $m\ge 0$. In this situation, one can do better, as shown by the following positive mass theorem with charge.
	\begin{theorem}\label{PMTC-Hyp}
		Let $(M^n,g,E)$, $n\geq 3$, be a complete time-symmetric spin initial data with charge containing at least one asymptotically hyperbolic end. Assume that the integrability condition (\ref{IntegrabilityHyperbolic}) holds, as well as the dominant energy condition (\ref{DEC-TS}) then $m\geq|Q|$.
	\end{theorem}
	We remark that, besides hyperbolic space, one might look for additional equality models within the Reissner--Nordström--AdS family. However, in the spherical case, the condition $m=|Q|$ does not give rise to a complete time-symmetric initial data set, but instead leads to a naked
	singularity, as shown in \cite{ChenLiZhou,WangZ}. Such configurations are therefore excluded from the class of admissible equality cases in our setting. From this perspective, introducing a compact inner boundary can be viewed as a way to discard the singular region and to investigate whether positivity can still be obtained in the time--symmetric setting. The general boundary results established above allow us to address this question.
	\begin{theorem}\label{PMTCB-Hyp-TS}
		Let $(M^n,g,E)$, $n\geq 3$, be a complete time--symmetric spin initial data set with charge containing at least one asymptotically hyperbolic end. Assume that
		the integrability condition \eqref{IntegrabilityHyperbolic} and the dominant energy condition \eqref{DEC-TS} hold, and that $(M^n,g)$ has a compact inner
		boundary $\partial M$. If one of the following boundary conditions is satisfied:
		\begin{enumerate}
			\item [{\rm (a)}]the mean curvature of $\partial M$ satisfies $H\le 0$;
			\item [{\rm (b)}] $\partial M$ has positive Yamabe invariant and satisfies
			\begin{align*}
				\frac{1}{n-1}\Big(H+(n-1)\left|E_\nu\right|\Big)\le
				\sqrt{\frac{1}{(n-1)(n-2)}\,\mathcal{Y}(\partial M,[\mathfrak{g}])\,
					{\rm Vol}(\partial M,\mathfrak{g})^{-\frac{2}{n-1}}+1},
			\end{align*}
		\end{enumerate}
		then the mass--charge inequality $m\ge |Q|$ holds.
	\end{theorem}
	Exterior coordinate spheres in the Reissner--Nordstr\"om--AdS manifold do not satisfy condition~{\rm (a)}, since their mean curvature is positive. 
	One might therefore hope that condition~{\rm (b)}, which allows for inner boundaries with positive mean curvature, could apply in this case. 
	This is however not the case, and the present result does not cover such exterior domains.
	
	We emphasize that no attempt is made here to analyze the equality case in the above results. This question will be addressed in a separate work \cite{Raulot17}.
	
	
	\section{The Schr\"odinger-Lichnerowicz formula for initial data sets with charge}\label{SL-Section}
	
	
	\subsection{Preliminaries on spinors}\label{IntroSpinors}
	
	Let $(M^n,g,K,E)$ be an initial data set with charge endowed with a spin structure. There exists a smooth Hermitian vector bundle over $M$, the spinor bundle, denoted by $\SB$, whose sections are called spinor fields. The Hermitian scalar product is denoted by $\<\,,\,\>$. Moreover, the tangent bundle $TM$ (in fact the Clifford bundle) acts on $\SB$ by Clifford multiplication $X\otimes \psi\mapsto c_g(X)\psi$ satisfying
	\begin{align*}
		c_g(X) c_g(Y)+c_g(Y)c_g(X)=-2g(X,Y)\,{\rm Id}
	\end{align*}
	and which is skew-Hermitian with respect to the Hermitian scalar product on $\SB$ that is 
	\begin{align*}
		\<c_g(X)\varphi,\psi\>  =  -\<\varphi,c_g(X)\psi\>
	\end{align*}
	for any tangent vector fields $X$ and any spinor fields $\varphi$, $\psi\in\Gamma(\SB)$. On the other hand, the Riemannian Levi-Civita connection $\nabla$ lifts to the so-called spin Levi-Civita connection, also denoted by $\nabla$, and defines a metric covariant derivative on $\SB$ that preserves the Clifford multiplication. This means that 
	\begin{align*}
		X\<\varphi,\psi\>  =  \<\nabla_X\varphi,\psi\>+\<\varphi,\nabla_X\psi\>
	\end{align*}
	and 
	\begin{align*}
		\nabla_X(c_g(Y)\varphi)  =  c_g\left(\nabla_XY\right)\varphi+c_g(Y)\nabla_X\varphi
	\end{align*}
	hold for all $X$, $Y\in\Gamma(TM)$ and $\varphi$, $\psi\in\Gamma(\SB)$. A quadruplet $(\SB,\<\,,\,\>,\nabla,\,\cdot\,)$ satisfying these properties is called a Dirac bundle. The Dirac operator is then the first order elliptic differential operator acting on $\SB$ locally defined by
	\begin{align*}
		D\varphi=\sum_{j=1}^nc_g(e_j)\nabla_{e_j}\varphi
	\end{align*}
	for $\varphi\in\Gamma(\SB)$. Here, and in all this work, $\{e_1,\cdots,e_n\}$ is a local orthonormal frame on $(M^n,g)$. 
	
	When $n$ is even, there exists a chirality operator, namely an endomorphism $\gamma$ of $\SB$ such that
	\begin{align}\label{ChiralityProperties}
		\<\gamma\varphi,\psi\>=\<\varphi,\gamma\psi\>,\quad\gamma^2={\rm Id},\quad \big\{ X,\gamma\big\}=0,\quad \nabla\gamma=0
	\end{align}
	where $\big\{ X,\gamma\big\}\varphi:=c_g(X)\gamma\varphi+\gamma\left(c_g(X)\varphi\right)$ for all $\varphi$, $\psi\in\Gamma(\SB)$ and $X\in\Gamma(TM)$. Such an operator is given by the Clifford action of the volume element of the spinor bundle $\SB$. When n is odd, a chirality operator in general does not exist. To overcome this, we consider $\SB\oplus\SB$, the direct sum of two copies of the spinor bundle over $M$, equipped with the Hermitian metric
	\begin{align*}
		\<(\varphi_1,\varphi_2),(\psi_1,\psi_2)\>:=\<\varphi_1,\psi_1\>+\<\varphi_2,\psi_2\>
	\end{align*}
	on which we define the Clifford action by
	\begin{align*}
		\mathfrak{c}_g(X)(\varphi_1,\varphi_2):=\left(c_g(X)\varphi_1,-c_g(X)\varphi_2\right)
	\end{align*}
	and the associated connection by $\nabla\oplus\nabla$. The endomorphism 
	\begin{align*}
		\gamma(\varphi_1,\varphi_2):=(\varphi_2,\varphi_1)
	\end{align*}
	then satisfies the same chirality properties (\ref{ChiralityProperties}) as in the even case. Thus by letting
	$$
	(\TSB,\<\,,\,\>,\nabla,\,c_g\,):=
	\left\lbrace
	\begin{array}{ll}
		(\SB,\<\,,\,\>,\nabla,\,c_g\,) & \text{ if } n \text{ is even}\\
		(\SB\oplus\SB,\<\,,\,\>\oplus\<\,,\,\>,\nabla\oplus\nabla,\,\mathfrak{c}_g\,) & \text{ if } n \text{ is odd}
	\end{array}
	\right.
	$$
	we obtain a Dirac bundle on which a chirality operator $\gamma\in\Gamma({\rm End\,}\TSB)$ satisfying (\ref{ChiralityProperties}) always exists. In the following, $D$ denotes the associated Dirac operator which satisfies the well-known Schr\"odinger-Lichnerowicz formula 
	\begin{align}\label{SL-Formula}
		D^2\varphi=\nabla^*\nabla\varphi+\frac{R_g}{4}\varphi
	\end{align}
	for all $\varphi\in\Gamma(\TSB)$. Here $\nabla^*$ is the $L^2$-formal adjoint of the connection $\nabla$ and the rough Laplacian is locally given by
	\begin{align}\label{RL-Local}
		\nabla^*\nabla\varphi=-\sum_{j=1}^n \nabla_{e_j}\nabla_{e_j}\varphi.
	\end{align}
	For a bounded domain $\Omega\subset M$ with smooth boundary, we have 
	\begin{align}\label{Dirac-IPP}
		\int_\Omega\<D\varphi,\psi\>d\mu=\int_\Omega\<\varphi,D\psi\>d\mu+\oint_{\partial\Omega}\<c_g(\nu)\varphi,\psi\>d\sigma
	\end{align}
	where $\nu$ is the outward unit normal to $\partial\Omega$ and $d\mu$ (resp. $d\sigma$) is the Riemannian volume form of $\Omega$ (resp. $\partial\Omega$) with respect to $g$. 
	
	When $g=b$, we indicate explicitly the dependence on $b$ in the notation ($\mathcal{S}_b$ for the spinor bundle, $\nabla^b$ the Riemannian or spin Levi-Civita connection, $\gamma_b$ the chirality operator,...).

	\subsection{The modified connection}\label{ModifiedConnectionSection}
	
	In this section, we let
	\begin{align*}
		\widetilde{\nabla}^\pm_X\varphi:=\nabla_X\varphi-\frac{1}{2}c_g(E) c_g(X)\gamma\varphi+\frac{n-3}{2}g(E,X)\gamma\varphi+\frac{1}{2}c_g\left(K(X)\right)\gamma\varphi\pm \frac{i}{2}c_g(X)\varphi
	\end{align*}
	for all $X\in \Gamma(TM)$ and $\varphi\in\Gamma(\TSB)$, so that the associated Dirac operators are locally given by
	\begin{align*}
		\widetilde{D}^\pm\varphi:=\sum_{j=1}^nc_g(e_j)\widetilde{\nabla}^\pm_{e_j}\varphi.
	\end{align*}
	
	To facilitate the computations, we re-express these connections in terms of the modified connection $\overline{\nabla}$ introduced in \cite{Raulot16} where positive energy theorems for spin initial data sets with charge are obtained. More precisely, we let 
	\begin{align}\label{ModifiedConnectionBar}
		\overline{\nabla}_X\varphi:=\nabla_X\varphi-\frac{1}{2}c_g(E)c_g(X)\gamma\varphi+\frac{n-3}{2}g(E,X)\gamma\varphi+\frac{1}{2}c_g\left(K(X)\right)\gamma\varphi
	\end{align}
	so that it is not difficult to check that the associated Dirac operator satisfies
	\begin{align}\label{RelationsDODBar}
		\overline{D}\varphi=D\varphi-\frac{1}{2}\left(c_g(E)+\Tr_gK\right)\gamma\varphi
	\end{align}
	for all $\varphi\in\Gamma(\TSB)$. Recall that this operator is an elliptic symmetric differential operator of order one and that it satisfies the following Schr\"odinger-Lichnerowicz formula 
	\begin{align}\label{ChargedSL}
		\overline{D}^2\varphi=\overline{\nabla}^*\overline{\nabla}\varphi+\mathcal{R}\varphi
	\end{align}
	for any $\varphi\in\Gamma(\TSB)$ and where $\mathcal{R}\in\Gamma({\rm End\,}\TSB)$ is defined by
	\begin{align*}
		\mathcal{R}\varphi=\frac{1}{2}\big(\mu\varphi+\varpi\gamma\varphi+c_g(J^\sharp)\gamma\varphi\big).
	\end{align*}
	It turns out that 
	\begin{align*}
		\widetilde{\nabla}^\pm_X\varphi=\overline{\nabla}_X\varphi\pm \frac{i}{2}c_g(X)\varphi
	\end{align*}
	so that 
	\begin{align}\label{DMT-DMB}
		\widetilde{D}^\pm\varphi=\overline{D}\varphi\mp i\frac{n}{2}\varphi
	\end{align}
	for all $X\in\Gamma(TM)$ and $\varphi\in\Gamma(\TSB)$.
	
	\begin{remark}\label{IPP-DM}
		When $\Omega$ is a bounded domain with smooth boundary $\partial\Omega$ in a spin initial data set with charge $(M^n,g,K,E)$, we deduce from (\ref{Dirac-IPP}) that 
		\begin{align*}
			\int_\Omega\<\widetilde{D}^\pm\varphi,\psi\>d\mu=\int_\Omega\<\varphi,\widetilde{D}^\mp\psi\>d\mu+\oint_{\partial\Omega}\<c_g(\nu)\varphi,\psi\>d\sigma
		\end{align*}
		for all $\varphi$, $\psi\in\Gamma(\TSB)$ and where $\nu$ is the outward unit normal to $\partial \Omega$. In particular, we have $\big(\widetilde{D}^\pm\big)^\ast=\widetilde{D}^\mp$ so that the operators $\widetilde{D}^\pm$ are not symmetric. 
	\end{remark}
	
	We can now state and prove the main formula of this section.
	\begin{theorem}\label{SLforEM}
		Let $(M^n,g,K,E)$ be a spin initial data set with charge, then
		\begin{align*}
			\big(\widetilde{D}^\pm\big)^\ast\widetilde{D}^\pm\varphi=\big(\widetilde{\nabla}^\pm\big)^\ast\widetilde{\nabla}^\pm\varphi+\widetilde{\mathcal{R}}\varphi
		\end{align*}
		for any $\varphi\in\Gamma(\TSB)$ and where $\widetilde{\mathcal{R}}\in\Gamma({\rm End\,}\TSB)$ is defined by
		\begin{align*}
			\widetilde{\mathcal{R}}\varphi=\frac{1}{2}\left(\left(\mu+\frac{n(n-1)}{2}\right)\varphi+\varpi\gamma\varphi+c_g(J^\sharp)\gamma\varphi\right).
		\end{align*}
	\end{theorem}

	{\it Proof.} We first note that it follows from (\ref{DMT-DMB}) and Remark \ref{IPP-DM} that 
	\begin{align*}
		\big(\widetilde{D}^+\big)^\ast\widetilde{D}^+\varphi=\widetilde{D}^-\widetilde{D}^+\varphi=\overline{D}^2\varphi+\frac{n^2}{4}\varphi.
	\end{align*}
	Combining this identity with the charged Schr\"odinger-Lichnerowicz formula (\ref{ChargedSL}) leads to 
	\begin{align}\label{FirstStepHSL}
		\big(\widetilde{D}^+\big)^\ast\widetilde{D}^+\varphi=\overline{\nabla}^*\overline{\nabla}\varphi+\left(\mathcal{R}+\frac{n^2}{4}\right)\varphi.
	\end{align}
	On the other hand, we compute for all $\varphi$, $\psi\in\Gamma(\TSB)$ that 
	\begin{eqnarray*}
		\<\widetilde{\nabla}^+\varphi,\widetilde{\nabla}^+\psi\> & = & \<\overline{\nabla}\varphi,\overline{\nabla}\psi\>+\frac{i}{2}\sum_{j=1}^n\<c_g(e_j)\varphi,\overline{\nabla}_{e_j}\psi\>+\frac{i}{2}\<\overline{D}\varphi,\psi\>+\frac{n}{4}\<\varphi,\psi\>\\
		& = & \<\overline{\nabla}^\ast\overline{\nabla}\varphi,\psi\>+\div(\xi_1)+\frac{i}{2}\sum_{j=1}^n\<c_g(e_j)\varphi,\overline{\nabla}_{e_j}\psi\>+\frac{i}{2}\<\overline{D}\varphi,\psi\>+\frac{n}{4}\<\varphi,\psi\>
	\end{eqnarray*}
	where $\xi_1\in\Gamma(TM)$ is the vector field on $M$ defined by $g(X,\xi_1)=\<\overline{\nabla}_X\varphi,\psi\>$ for all $X\in\Gamma(TM)$. The second equality is obtained from the first by virtue of \cite[Remark 3]{Raulot16}. Now using (\ref{ModifiedConnectionBar}) we write
	\begin{eqnarray*}
		\sum_{j=1}^n\<c_g(e_j)\varphi,\overline{\nabla}_{e_j}\psi\> & =  & \underbrace{\sum_{j=1}^n\<c_g(e_j)\varphi,\nabla_{e_j}\psi\>}_{(1)}-\underbrace{\frac{1}{2}\sum_{j=1}^n\<c_g(e_j)\varphi,c_g(E)c_g(e_j)\gamma\psi\>}_{(2)}\\& +& \underbrace{\frac{n-3}{2}\sum_{j=1}^n\<c_g(e_j)\varphi,E^j\gamma\psi\>}_{(3)}+\underbrace{\frac{1}{2}\sum_{j=1}^n\<c_g(e_j)\varphi,K(e_j)\gamma\psi\>}_{(4)}
	\end{eqnarray*}
	where $E^j:=g(E,e_j)$ for $j=1,\cdots,n$. A straightforward computation of the first term yields
	\begin{align*}
		(1) = \div_g(\xi_2)-\<D\varphi,\psi\>
	\end{align*}
	where $\xi_2\in\Gamma(TM)$ is the vector field on $M$ defined by $g(\xi_2,X)=\<c_g(X)\varphi,\psi\>$ for all $X\in\Gamma(TM)$. 
	On the other hand, the Clifford rule ensures that 
	\begin{align*}
		(2)  =  \frac{1}{2}\sum_{j=1}^n\<c_g(e_j)c_g(E)c_g(e_j)\varphi,\gamma\psi\> =  -\frac{1}{2}\sum_{j=1}^n\<c_g(E)c_g(e_j) c_g(e_j)\varphi,\gamma\psi\>-\<c_g(E)\varphi,\gamma\psi\>
	\end{align*}
	and so $(2)=\frac{n-2}{2}\<c_g(E)\varphi,\gamma\psi\>$. It is also easy to see that $(3)=\frac{n-3}{2}\<c_g(E)\varphi,\gamma\psi\>$. Finally, we also see that 
	\begin{align*}
		(4)=-\frac{1}{2}\sum_{j=1}^n\<c_g\left(K(e_j)\right)c_g(e_j)\varphi,\gamma\psi\>= \frac{1}{2}\left(\Tr_gK\right)\<\varphi,\gamma\psi\>
	\end{align*}
	since $K$ is symmetric. Putting together the expressions of these four terms and using equality (\ref{RelationsDODBar}) we deduce that
	\begin{align*}
		\sum_{j=1}^n\<c_g(e_j)\varphi,\overline{\nabla}_{e_j}\psi\> = \div_g(\xi_2)-\<\overline{D}\varphi,\psi\>
	\end{align*}
	and so 
	\begin{align*} 
		\<\widetilde{\nabla}^+\varphi,\widetilde{\nabla}^+\psi\>  = \div_g(\xi)+\<\overline{\nabla}^\ast\overline{\nabla}\varphi+\frac{n}{4}\varphi,\psi\>
	\end{align*}
	where $\xi\in\Gamma(TM)$ is the vector field defined by $g(\xi,X)=\<\widetilde{\nabla}^+_X\varphi,\psi\>$ for all $X\in\Gamma(TM)$. We have thus shown that 
	\begin{align*}
		\big(\widetilde{\nabla}^+\big)^\ast\widetilde{\nabla}^+\varphi=\overline{\nabla}^\ast\overline{\nabla}\varphi+\frac{n}{4}\varphi
	\end{align*}
	which, when combined with (\ref{FirstStepHSL}), leads to the claimed formula.
	\qed
	
	\begin{remark}\label{DEC-SL}
		The curvature endomorphism $\widetilde{\mathcal{R}}$ satisfies
		\begin{align*}
			\<\widetilde{\mathcal{R}}\varphi,\varphi\>
			\geq \frac{1}{2}\left(\left(\mu+\frac{n(n-1)}{2}\right)
			-\sqrt{|J|_g^2+\varpi^2}\right)|\varphi|^2
		\end{align*}
		for all $\varphi\in\Gamma(\TSB)$. Under the dominant energy condition~\eqref{DEC-General}, this implies that $\widetilde{\mathcal{R}}$ is non-negative.
	\end{remark}
	
	\begin{remark}\label{IPP-Nablapm}
		In the preceding proof, we proved the pointwise equality:
		\begin{align*}
			\<\big(\widetilde{\nabla}^\pm\big)^\ast\widetilde{\nabla}^\pm\varphi,\psi\>= \<\widetilde{\nabla}^\pm\varphi,\widetilde{\nabla}^\pm\psi\>-\div_g(\xi^\pm)
		\end{align*}
		where $\xi^\pm\in\Gamma(TM)$ is the vector field defined by $g(\xi^\pm,X)=\<\widetilde{\nabla}^\pm_X\varphi,\psi\>$ for all $X\in\Gamma(TM)$ and where $\varphi$, $\psi\in\Gamma(\TSB)$. Then if $\Omega$ is a bounded domain with smooth boundary in a complete spin initial data set $(M^n,g,K,E)$, one can apply the divergence formula to obtain the following integration by parts formula:
		\begin{align*}
			\int_\Omega\<\big(\widetilde{\nabla}^\pm\big)^\ast\widetilde{\nabla}^\pm\varphi,\psi\>d\mu=\int_\Omega\<\widetilde{\nabla}^\pm\varphi,\widetilde{\nabla}^\pm\psi\>d\mu-\oint_{\partial\Omega}\<\widetilde{\nabla}^\pm_\nu\varphi,\psi\>d\sigma.
		\end{align*}
	\end{remark}
	
	For compact domains, we obtain the general integral formula: 
	\begin{corollary}\label{ReillyHyperbolic}
		Let $\Omega$ be a compact domain with boundary in a complete spin initial data with charge $(M^n,g,K,E)$. Then, for any $\varphi\in\Gamma(\TSB)$, we have
		\begin{align*}
			\int_\Omega\left(|\widetilde{\nabla}^\pm\varphi|^2+\<\widetilde{\mathcal{R}}\varphi,\varphi\>-|\widetilde{D}^\pm\varphi|^2\right)d\mu=\oint_{\partial\Omega}\<\widetilde{L}^\pm_\nu\varphi,\varphi\>d\sigma.
		\end{align*}
		where $\nu$ is the outward unit normal to $\partial\Omega$ and $\widetilde{L}^\pm_\nu\varphi:=\widetilde{\nabla}^\pm_\nu\varphi+c_g(\nu)\widetilde{D}^\pm\varphi$.
	\end{corollary}	
	
	
	\section{Proof of the main results}\label{PMR}
	
	
	In this part, we begin by recalling the construction of extrinsic imaginary Killing spinors on hyperbolic space, essentially going back to \cite{Baum2}. When $n$ is even, these coincide with the classical imaginary Killing spinors; when $n$ is odd, one must take into account the doubling construction described in Section~\ref{IntroSpinors}. For convenience, set
	\begin{align*}
		\Sigma_n :=
		\begin{cases}
			\mathbb{C}^{2^m} & \text{if } n=2m,\\[0.2cm]
			\mathbb{C}^{2^m}\oplus\mathbb{C}^{2^m} & \text{if } n=2m+1.
		\end{cases}
	\end{align*}
	Using the unit ball model $(\mathbb{B}^n,\omega^{-2}\delta)$ of $\mathbb{H}^n$, where $\omega(x)=(1-|x|^2)/2$, we identify the trivial spinor bundle $\mathbf{S}_\delta = \mathbb{B}^n\times\mathbb{C}^{[n/2]}$ with the spinor bundle $\mathbf{S}_b$ of the hyperbolic metric $b$. For each constant spinor $u\in\Gamma(\mathbf{S}_\delta)$, we define
	\begin{align*}
		\varphi_u^\pm(x)
		= \omega(x)^{-1/2}\,\big(1 \mp i\, c_\delta(x)\big)\,u,
	\end{align*}
	which gives an imaginary Killing spinor on $(\mathbb{H}^n,b)$. Here we omit the identification between $\mathbf{S}_\delta$ and $\mathbf{S}_b$. 
	
	When $n=2m$, the map $u\mapsto\zeta_u:=\varphi^+_u$ provides a complete parametrisation of the space $\EH^+$ of (extrinsic) imaginary Killing spinors.  When $n=2m+1$, we use the doubled bundle $\mathcal{S}_b=\mathbf{S}_b\oplus\mathbf{S}_b$ and the chirality operator $\gamma_b$ of Section~\ref{IntroSpinors} to define, for $u=(v,w)\in\Sigma_n$,
	\begin{align*}
		\zeta_u := 
		\begin{pmatrix}
			\varphi_v^+ \\[0.1cm]
			\varphi_w^- 
		\end{pmatrix}
		\in \EH^+.
	\end{align*}
	Thus, in all dimensions, the assignment
	\begin{align*}
		\zeta : u\in\Sigma_n \longmapsto \zeta_u\in \EH^+
	\end{align*}
	is well-defined. From each such $u\in\Sigma_n$ we form the triple $(V_u,\alpha_u,f_u)$ defined by
	\begin{align*}
		\begin{cases}
			V_u = |\zeta_u|^2_b,\\[0.1cm]
			\alpha_u(X) = \langle c_b(X)\gamma_b\zeta_u,\zeta_u\rangle_b,\\[0.1cm]
			f_u = \langle\gamma_b\zeta_u,\zeta_u\rangle_b,
		\end{cases}
		\qquad X\in\Gamma(T\mathbb{H}^n),
	\end{align*} 
	which is precisely the quadratic map $\mathcal{K}$ appearing in Theorem~\ref{PET-Spherical}.
	\begin{proposition}
		The  quadratic map $\mathcal{K}$ takes values in the cone ${\mathscr C}_b$.
	\end{proposition}
	
	{\it Proof.} The facts that $V_u\in\mathcal{N}_b$ and $\alpha_u\in\mathfrak{Kill}(\mathbb{H}^n)$ are well-known (see for example \cite{Maerten2,ChruscielMaertenTod}) and can be straightforwardly deduced from the Killing equation (\ref{EKS}). For the last point, we compute using (\ref{EKS}) that 
	\begin{align*}
		X(f_u)= 2{\rm Re}\<\nabla^b_X(\gamma_b\zeta_u),\zeta_u\>=2{\rm Re}\<\gamma_b\left(\nabla^b_X\zeta_u\right),\zeta_u\>=-{\rm Re}\left(i\<\gamma_b\left(c_b(X)\zeta_u\right),\zeta_u\>\right)
	\end{align*}
	for all $X\in\Gamma(T\mathbb{H}^n)$. On the other hand, it follows from the properties (\ref{ChiralityProperties}) that
	\begin{align*}
		\overline{\<\gamma_b\left(c_b(X)\zeta_u\right),\zeta_u\>}=\<\gamma_b\left(c_b(X)\zeta_u\right),\zeta_u\>
	\end{align*}
	and so $X(f_u)=0$ which implies that $f_u$ is a constant. 
	\qed
	
	\medskip
	
	Now we remark that the boundary contribution of the integrated Schr\"o\-ding\-er-Li\-chne\-ro\-wicz formula can be identified with the global charges $\Xi(V,\alpha,f)$ for all $(V,\alpha,f)\in\mathscr{C}_b$. Indeed, since $(M^n,g)$ is asymptotically AdS, we can identify the spinor bundle over $(R,\infty)\times\mathbb{S}^{n-1}$ with the one over $M_{ext}$ by using the chart at infinity. We will denote by $\mathcal{A}$ this identification. So, through this identification, every extrinsic imaginary Killing spinor $\zeta\in\Gamma(\EH^\pm)$ is mapped to a spinor field over $M_{ext}$, denoted by $\mathcal{A}\zeta$. We first remark the following property:
	\begin{proposition}\label{CEM-Spinor}
		For all $\zeta\in\Gamma(\EH^\pm)$, we have 
		\begin{align*}
			\lim_{r\rightarrow\infty}\oint_{S_r}\<\widetilde{L}^\pm_{\nu_r}\phi,\phi\>d\sigma_r=\frac{n-1}{2}\omega_{n-1}\Xi\left(V_{\zeta},\alpha_{\zeta},f_{\zeta}\right)
		\end{align*}
		where $\phi:=\eta\mathcal{A}\zeta\in\Gamma(\TSB)$ and $\eta$ is a cut-off function that vanishes outside of $M_{ext}$ and is equal to $1$ for large $r$.
	\end{proposition}
	{\it Proof.} First we write the boundary operator $\widetilde{L}^\pm_{\nu_r}$ in terms of the corresponding operator
	\begin{align*}
		\widehat{L}^{\pm}_{\nu_r}\phi:=\widehat{\nabla}^\pm_{\nu_r}\phi+c_g(\nu_r)\widehat{D}^\pm\phi
	\end{align*}
	where the connection $\widehat{\nabla}^\pm$ is defined by
	\begin{align*}
		\widehat{\nabla}^\pm_X\phi:=\nabla_X\phi+\frac{1}{2}c_g\left(K(X)\right)\gamma\phi\pm\frac{i}{2}c_g(X)\phi
	\end{align*}
	and whose associated Dirac operator is denoted by $\widehat{D}^\pm$. Note that a straightforward computation gives 
	\begin{align*}
		\widetilde{D}^\pm\phi=\widehat{D}^\pm\phi-\frac{1}{2}c_g(E)\gamma\phi
	\end{align*}
	leading to the relation 
	\begin{align*}
		\widetilde{L}^{\pm}_{\nu_r}\phi  =  \widehat{L}^{\pm}_{\nu_r}\phi+\frac{n-1}{2}E_{\nu_r}\gamma\phi.
	\end{align*}
	From \cite[Proposition 3.2]{Maerten2} we get that
	\begin{align*}
		\lim_{r\rightarrow\infty}\oint_{S_r}\<\widehat{L}^\pm_{\nu_r}\phi,\phi\>d\sigma_r=\frac{1}{4}\lim_{r\rightarrow\infty}\oint_{S_r}\Big(\mathbb{U}_1(V_\zeta)+\mathbb{U}_2(\alpha_\zeta)\Big)\left(\nu_{b_r}\right)d\sigma_{b_r}.
	\end{align*}
	On the other hand, since $\mathcal{A}$ is a pointwise isometry between the spinor bundle over the end with $\gamma\circ\mathcal{A}=\mathcal{A}\circ\gamma_b$, a direct computation using the asymptotic behavior of the metric $g$ gives
	\begin{align*}
		\frac{n-1}{2}\lim_{r\rightarrow\infty}\oint_{S_r}E_{\nu_r}\<\gamma\phi,\phi\>d\sigma_r=\frac{n-1}{2}\omega_{n-1}\lim_{r\rightarrow\infty}\oint_{S_r}\mathbb{U}_3(f_\zeta)\left(\nu_{b_r}\right)d\sigma_{b_r}
	\end{align*}
	and the claim of the proposition follows.
	\qed
	
	The second step of the proof is to show that the modified Dirac operator $\widetilde{D}^\pm$ is an isomorphism between a suitable Hilbert space and
	$L^2$, the space of square-integrable sections of $\TSB$. For this, we adopt the framework developed by Bartnik and Chru\'sciel \cite{BartnikChrusciel}.  We first observe that under our assumptions, we have a {\it weighted Poincar\'e inequality} for the connections $\widetilde{\nabla}^\pm$. More precisely, this means that there exists $w\in L^1_{loc}$ with ${\rm ess\,inf}_\Omega w>0$ for all relatively compact $\Omega$ in $M$ such that for all $\varphi\in C^1_c$, the space of compactly supported $C^1$ spinor fields on $M$, we have 
	\begin{align}\label{wPi}
		\int_M|\varphi|^2w\,d\mu\leq\int_M |\widetilde{\nabla}^\pm\varphi|^2d\mu.
	\end{align}
	Indeed, since the symmetric part $\widetilde{\Gamma}^\pm_S$ of the modified connection
	$\widetilde{\nabla}^\pm$ is given by
	\begin{align*}
		\widetilde{\Gamma}^\pm_S(X)=-\frac{1}{2}\big(c_g\big(K(X)\big)+(n-2)g(E,X)\big)\gamma\mp\frac{i}{2}c_g(X)\in\Gamma\big({\rm End\,}\TSB\big)
	\end{align*}
	it follows from the decay assumptions \eqref{AsymptoticDecay} that for $n\geq 3$ the hypotheses of \cite[Theorem~9.10]{BartnikChrusciel} are satisfied.
	
	
	\subsection{The boundaryless case: proof of Theorem~\ref{PET-Spherical}}\label{BoundaryLessCase}
	
	In this situation, we have the following {\it Schr\"odinger-Lichnerowicz estimate} for the Dirac operator $\widetilde{D}^+$ (in the sense of \cite{BartnikChrusciel}) 
	\begin{align}\label{SL-estimate}
		\int_M|\widetilde{\nabla}^+\phi|^2d\mu\leq \int_M|\widetilde{D}^+\phi|^2d\mu
	\end{align}
	for all $\phi\in C^1_c$. This follows directly by integrating the formula of Theorem \ref{SLforEM} and using the dominant energy condition (see Remark \ref{DEC-SL}). 
	This last property also implies that 
	\begin{align*}
		||\phi||^2 :=\int_{M}\Big(|\widetilde{\nabla}^+\phi|^2+\<\widetilde{\mathcal{R}}\phi,\phi\>\Big)d\mu
	\end{align*}
	defines a norm on $C^1_c$. Therefore the space
	\begin{align*}
		\mathbf{H}:=||\,.\,||-{\rm completion\,\,of\,\,}C^1_c
	\end{align*}
	is a Hilbert space. The Poincar\'e inequality (\ref{wPi}) ensures that $\mathbf{H}$ embeds continuously in $H^1_{loc}$. In particular, it implies that any $\phi\in\mathbf{H}$ can be represented by a spinor field in $H^1_{loc}$. Consider now the sesquilinear form defined by
	\begin{align*}
		\alpha(\phi,\psi):=\int_M\<\widetilde{D}^+\phi,\widetilde{D}^+\psi\>d\mu
	\end{align*}
	for $\phi$, $\psi\in\mathbf{H}$. From Lemma $8.5$ in \cite{BartnikChrusciel}, we get that the map $\phi\in\mathbf{H}\mapsto \widetilde{D}^+\phi\in L^2$ is continuous and so $\alpha$ is also continuous on $\mathbf{H}\times\mathbf{H}$. Moreover, using the weighted Poincar\'e inequality (\ref{wPi}) and the Schr\"odinger-Lichnerowicz estimate (\ref{SL-estimate}), we immediately see that $\alpha$ is coercive on $\mathbf{H}$. For $\chi\in L^2$, consider the continuous linear form
	\begin{align*}
		F_\chi(\phi)=\int_M\<\chi,\widetilde{D}^+\phi\>d\mu
	\end{align*}
	on $\mathbf{H}$. The Lax-Milgram theorem then implies the existence of a unique $\xi\in\mathbf{H}$ such that $F_\chi(\phi)=\alpha(\xi,\phi)$ for all $\phi\in\mathbf{H}$. Equivalently, $\xi_0:=\widetilde{D}^+\xi-\chi\in L^2$ is a weak solution of $\widetilde{D}^-\xi_0=0$ since we know from Remark \ref{IPP-DM} that $(\widetilde{D}^+)^\ast=\widetilde{D}^-$. From the ellipticity of $\widetilde{D}^-$, we conclude that $\xi_0\in\mathbf{H}\cap L^2$ is in fact a strong solution of this equation. On the other hand, since it holds that 
	\begin{align}\label{UnicityStep}
		\widetilde{D}^-\xi_0=\widetilde{D}^+\xi_0+in\xi_0
	\end{align}
	we get that $(\widetilde{D}^+)^k\xi_0\in L^2$ for all $k\in\mathbb{N}$. Therefore, this shows that 
	\begin{align*}
		\int_M|\widetilde{D}^+\xi_0|^2d\mu=\int_M\<\widetilde{D}^-\widetilde{D}^+\xi_0,\xi_0\>d\mu=\int_M\<\widetilde{D}^+\widetilde{D}^-\xi_0,\xi_0\>d\mu=0
	\end{align*}
	and so $\widetilde{D}^+\xi_0=0$ which, with (\ref{UnicityStep}), ensures that $\xi_0=0$ and so $\xi\in\mathbf{H}$ is the unique solution of $\widetilde{D}^+\xi=\chi$ for $\chi\in L^2$. To summarize, we proved the following result:
	\begin{proposition}\label{DiracIsomorphismWB}
		Under the assumptions of Theorem \ref{PET-Spherical}, the operator $\widetilde{D}^+:\mathbf{H}\rightarrow L^2$ is an isomorphism.
	\end{proposition}
	Now we can apply the classical Witten argument to conclude. Take $\zeta\in\Gamma(\EH^+)$ and consider $\phi:=\eta\mathcal{A}\zeta\in\Gamma(\TSB)$ as in Proposition \ref{CEM-Spinor}. It follows from our asymptotic assumptions (\ref{AsymptoticDecay}) that $\widetilde{D}^+\phi$ is $L^2$ on $(M^n,g)$ so that Proposition \ref{DiracIsomorphismWB} ensures the existence of a unique $\xi\in\mathbf{H}$ such that $\widetilde{D}^+\xi=-\widetilde{D}^+\phi$. In other words, the spinor field $\Theta:=\xi+\phi$ is $\widetilde{D}^+$-harmonic. If $\xi\in C^1_c$, Proposition~\ref{CEM-Spinor} would directly yield
	\begin{eqnarray}\label{mass-formula}
		\frac{n-1}{2}\omega_{n-1}\Xi\left(V_\zeta,\alpha_\zeta,f_\zeta\right) & = & \lim_{r\rightarrow\infty}\oint_{S_r}\<\widetilde{L}^+_{\nu_r}\Theta,\Theta\>d\sigma_r\nonumber\\
		& = & \int_{M}\left(|\widetilde{\nabla}^+\Theta|^2+\<\widetilde{\mathcal{R}}\Theta,\Theta\>\right)d\mu\geq 0
	\end{eqnarray}
	because of Remark \ref{DEC-SL}. Actually, one can then show that the previous equality holds for $\xi\in\mathbf{H}$ since the right-hand side of (\ref{mass-formula}) is continuous on $\mathbf{H}$ (and $C^1_c$ is dense in $\mathbf{H}$). We finally get that the linear form $\Xi$ is non-negative on the cone $\mathscr{C}_b$ and this concludes the proof of Theorem \ref{PET-Spherical}.
	
	\subsection{The non-empty boundary case: proof of Theorem~\ref{PositiveETS-Boundary}}\label{BoundaryPart1}
	
	From now on we assume that $M$ has a compact inner boundary. To control the additional boundary term coming from this boundary, we need to introduce some notations. On the restricted spinor bundle $\ETSB_g:=\mathcal{S}_{g|\partial M}$ over $\partial M$, we define for $X\in\Gamma(T\partial M)$ the linear connection 
	\begin{align*}
		\nb_X\varphi :=\nabla_X\varphi+\frac{1}{2} c_g\left(A(X)\right)c_g(\nu)\varphi
	\end{align*}
	where $A(X):=\nabla_X\nu$ is the Weingarten map of $\partial M$ in $M$ and the associated Dirac operator, denoted by $\D$, is then as usual locally given by
	\begin{align*}
		\D\varphi=\sum_{i=1}^{n-1}c_g(e_i)c_g(\nu)\nb_{e_i}\varphi.
	\end{align*}
	Here $\nu$ denotes the unit normal to $\partial M$ pointing toward infinity. Moreover, $H$ denotes the mean curvature of $\partial M$, $\Tr_{\partial M}(K)$ is the trace along $\partial M$ of $K$ and $K(\nu)^\top$ is the tangent part of the vector field $K(\nu)$ defined along $\partial M$ by 
	\begin{align*}
		g\big(K(\nu),X\big)=K(X,\nu)
	\end{align*}
	for all $X\in\Gamma(TM_{|\partial M})$. We then prove:
	\begin{proposition}\label{BoundaryExpression}
		For all $\varphi\in\Gamma(\ETSB_{g})$, the following identity holds
		\begin{align*}
			\widetilde{L}^\pm_\nu\varphi=-\D^\pm\varphi-\frac{1}{2}\mathcal{H}\varphi
		\end{align*}
		where $\D^\pm$ is the first order elliptic differential operator given by
		\begin{align*}
			\D^\pm\varphi=\D\varphi\pm\frac{n-1}{2}ic_g(\nu)\varphi
		\end{align*}
		and $\mathcal{H}\in\Gamma({\rm End\,}\ETSB_{g})$ is defined by
		\begin{align*}
			\mathcal{H}\varphi:=H\varphi+\Tr_{\partial M}(K)c_g(\nu)\gamma\varphi-c_g\big(K(\nu)^\top\big)\gamma\varphi-(n-1)E_\nu\gamma\varphi.
		\end{align*}
	\end{proposition}
	
	{\it Proof.} It suffices to observe that
	\begin{align*}
		\widetilde{L}^\pm_\nu\varphi=\overline{L}_\nu\varphi\mp\frac{n-1}{2}ic_g(\nu)\varphi
	\end{align*}
	where $\overline{L}_\nu\varphi:=\overline{\nabla}_\nu\varphi+c_g(\nu)\overline{D}\varphi$ and to use \cite[Proposition 3]{Raulot16} to conclude.
	\qed
	
	Using Corollary \ref{ReillyHyperbolic} on the compact domain $M_r$ whose boundary is the union of $\partial M$ and $S_r$, we deduce that:
	\begin{equation}\label{IntegralVersionBoundary}
		\begin{aligned}
			\int_{S_r} \< \widetilde{L}^\pm_{\nu_r} \varphi, \varphi \> \, d\sigma_r
			&= \int_{M_r} \Big( |\widetilde{\nabla}^\pm \varphi|^2
			+ \< \widetilde{\mathcal{R}} \varphi, \varphi \>
			- |\widetilde{D}^\pm \varphi|^2 \Big) \, d\mu \\
			&\quad - \oint_{\partial M} \< \D^\pm \varphi
			+ \tfrac{1}{2} \mathcal{H} \varphi, \varphi \> \, d\sigma.
		\end{aligned}
	\end{equation}
	for all $\varphi\in\Gamma(\TSB)$. To proceed as in the boundaryless case, one needs to impose boundary conditions to the operators $\widetilde{D}^\pm$ to ensure sufficiently good analytic properties. For this purpose, we consider the pointwise projections
	\begin{align*}
		\pi_\pm \varphi := \frac{1}{2} (\varphi\pm c_g(\nu) \gamma \varphi), \qquad \varphi \in \Gamma(\ETSB_g),
	\end{align*}
	which are well-known to define elliptic boundary conditions for the standard Dirac-Witten operator $D$ (see, for example, \cite{BartnikChrusciel}). Since ellipticity depends only on the principal symbol and $D$ and $\widetilde{D}^\pm$ share the same principal symbol, these conditions are also elliptic for $\widetilde{D}^\pm$. Now let 
	\begin{align*}
		\widetilde{C}^1_c:=\big\{\varphi\in C^1_c\,\big|\, \pi_-\varphi_{|\partial M_+}=0\text{ and }\pi_+\varphi_{|\partial M_-}=0\big\}
	\end{align*}
	where $\partial M = \partial M_+ \cup \partial M_-$, with $\partial M_\pm$ denoting the portion of the boundary on which $\theta_\pm \le 0$. Here $\theta_\pm$ denote the null expansions of $\partial M$ defined by (\ref{FuturePastTrapped}). Then we prove:
	\begin{proposition}
		Under the assumptions of Theorem \ref{PositiveETS-Boundary}, the Schr\"odinger-Lichnerowicz estimate (\ref{SL-estimate}) holds for all $\varphi\in\widetilde{C}_c^1$. 
	\end{proposition}
	{\it Proof.} A direct consequence of formula (\ref{IntegralVersionBoundary}) and the dominant energy condition is that 
	\begin{align*}
		\int_M|\widetilde{D}^+\varphi|^2d\mu\geq \int_{M}|\widetilde{\nabla}^+\varphi|^2d\mu - \oint_{\partial M}\<\D^+\varphi+\frac{1}{2}\mathcal{H}\varphi,\varphi\>d\sigma
	\end{align*}
	for all $\varphi\in\widetilde{C}^1_c$. Now it is not difficult to check using the definition of $\D$ and (\ref{ChiralityProperties}) that
	\begin{align*}
		\<\D\varphi,\varphi\> = \<\D(\pi_+\varphi),\pi_-\varphi\>+\<\D(\pi_-\varphi),\pi_+\varphi\>
	\end{align*}
	for all $\varphi\in\Gamma(\TSB)$ so that $\<\D\varphi,\varphi\>=0$ as soon as $\varphi\in\widetilde{C}^1_c$. One computes next that
	\begin{align*}
		i\< c_g(\nu)\varphi,\varphi\>
		= -2\,\mathrm{Im}\,\langle c_g(\nu)(\pi_+\varphi),\pi_-\varphi\rangle ,
	\end{align*}
	which in particular also vanishes identically for all $\varphi\in\widetilde{C}^1_c$. Similarly, it is straightforward to check that
	\begin{align*}
		\<c_g(\nu)\gamma\varphi,\varphi\>  =  |\pi_+\varphi|^2-|\pi_-\varphi|^2,\quad\<\gamma\varphi,\varphi\> =  2\,{\rm Re} \<\gamma(\pi_+\varphi),\pi_-\varphi\>
	\end{align*}
	and 
	\begin{align*} 
		\<c_g(X)\gamma\varphi,\varphi\>  = 2\,{\rm Re} \<c_g(X)\gamma(\pi_+\varphi),\pi_-\varphi\>
	\end{align*}
	for all $X\in\Gamma(T\partial M)$ and $\varphi\in\Gamma(\TSB)$. Thus it holds on $\partial M_\pm$ that $\<\mathcal{H}\varphi,\varphi\>=\theta_\pm|\pi_{\pm}\varphi|^2$ for all $\varphi\in\widetilde{C}^1_c$. Combining the previous estimates leads to 
	\begin{align*}
		\int_M|\widetilde{D}^+\varphi|^2d\mu\geq \int_{M}|\widetilde{\nabla}^+\varphi|^2d\mu -\frac{1}{2}\oint_{\partial M_+}\theta_+|\pi_+\varphi|^2d\sigma-\frac{1}{2}\oint_{\partial M_-}\theta_-|\pi_-\varphi|^2d\sigma
	\end{align*}
	and this concludes the proof since $\theta_\pm\leq 0$ on $\partial M_{\pm}$. 
	\qed
	
	By taking $\mathbf{H}$ to be the $||\cdot||$--completion of $\widetilde{C}^1_c$, the Green formula of Remark~\ref{IPP-DM} shows that the adjoint of $\widetilde{D}^+:\mathbf{H}\to L^2$ is $\widetilde{D}^-:\mathbf{H}\to L^2$. The proof of Theorem~\ref{PositiveETS-Boundary} then proceeds exactly as in the boundaryless case.
	
	\subsection{The non-empty boundary case: proof of Theorem~\ref{PositiveETS-Boundary-APS}}
	
	The proof of Theorem~\ref{PositiveETS-Boundary-APS} relies on Atiyah--Patodi--Singer type boundary conditions $\mathcal P_{+}^\pm$, following an approach introduced by Herzlich \cite{Herzlich1,Herzlich2} in the context of positive mass theorems. We briefly recall these conditions below and refer to \cite[Sections~3.2 and~3.3]{HijaziMontielRaulot4} for further details. As before, the key point is to establish a Schr\"odinger--Lichnerowicz estimate for the pair $(\widetilde{D}^\pm,\mathcal P_{+}^\pm)$.
	
	In Proposition~\ref{BoundaryExpression}, we introduced the operators $\D^\pm$ on the compact boundary $\partial M$. They are first-order elliptic and self-adjoint, with discrete, unbounded spectra of non-zero real eigenvalues, symmetric with respect to zero; in fact, the two spectra coincide. In particular, they satisfy
	\begin{align}\label{RelationDpm}
		\D^+\big(c_g(\nu)\psi\big) = -c_g(\nu)\D^-\psi
	\end{align}
	for all $\psi\in\Gamma(\ETSB_g)$. We denote by $(\mu_k)_{k\in\mathbb Z^\ast}$ their common spectrum, ordered so that $\mu_k>0$ for $k>0$ and $\mu_{-k}=-\mu_k$. We let $\mathcal P_{+}^\pm$ for the $L^2$-orthogonal projection onto the subspace spanned by the $\D^{\pm}$-eigenspinors corresponding to positive eigenvalues. As shown in \cite{BartnikChrusciel}, these projections define global elliptic boundary conditions for $\widetilde{D}^\pm$. Using the integration by parts formula of Remark~\ref{IPP-DM} together with~\eqref{RelationDpm}, one checks that the adjoint of $\widetilde{D}^\pm$ endowed with $\mathcal P_{+}^\pm$ is $\widetilde{D}^\mp$ endowed with $\mathcal P_{+}^\mp$.
	
	A further key point in the proof of Theorem~\ref{PositiveETS-Boundary-APS} is the control of the sign of the boundary contribution in~\eqref{IntegralVersionBoundary}. This requires a lower bound on the first positive eigenvalue $\mu_1$ of $\D^\pm$. Such a bound is obtained by relating $\mu_1$ to the first positive eigenvalue of the Dirac operator $\D$ on $\partial M$, for which several sharp estimates are available. In particular, applying the results of B\"ar \cite{Bar3} and Hijazi \cite{Hijazi2,Hijazi1} yields the inequality:
	\begin{align}\label{BarHijaziRaulot}
		\mu_1^2
		\ge
		\frac{n-1}{4(n-2)}\frac{\mathcal{Y}(\partial M,[\mathfrak{g}])}{{\rm Vol}(\partial M,\mathfrak{g})^{\frac{2}{n-1}}}
		+
		\frac{(n-1)^2}{4}
	\end{align}
	when the Yamabe invariant $\mathcal{Y}(\partial M,[\mathfrak{g}])$ of $(\partial M,\mathfrak{g})$ is positive. 
	
	We now introduce the space
	\begin{align*}
		C_c^{1,\pm}:=\big\{\varphi\in C^1_c\,\big|\, \mathcal P_{+}^\pm(\varphi_{|\partial M})=0\big\}.
	\end{align*}
	For any $\varphi\in C_c^{1,\pm}$, its boundary trace admits the spectral decomposition
	\begin{align*}
		\varphi_{|\partial M} = \sum_{k<0} A^\pm_k\psi_k^\pm,
		\quad\text{where}\quad
		A^\pm_k=\oint_{\partial M}\<\varphi,\psi_k^\pm\>d\sigma\in\mathbb{C}.
	\end{align*}
	Here $(\psi_k^\pm)_{k\in\mathbb{Z}^*}$ is an $L^2$-orthonormal basis of $\D^\pm$-eigenspinors. Using this decomposition, the boundary term can be estimated as:
	\begin{eqnarray*}
		\oint_{\partial M}
		\<
		\D^\pm\varphi+\frac{1}{2}\mathcal H_{\max}\varphi,\varphi
		\> d\sigma
		& \le &  
		\left(\frac{1}{2}\sup_{\partial M}\left(\mathcal H_{\max}\right)-\mu_1\right)
		\oint_{\partial M}|\varphi|^2d\sigma \leq 0.
	\end{eqnarray*}
    The last inequality follows from the lower bound~\eqref{BarHijaziRaulot} and the assumption~\eqref{YamabeHmax} in Theorem~\ref{PositiveETS-Boundary-APS}. Together with Theorem~\ref{SLforEM} and the dominant energy condition~\eqref{DEC-General}, this yields a Schr\"odinger--Lichnerowicz estimate for the pair $(\widetilde{D}^\pm,\mathcal{P}_{+}^\pm)$.
	
	The proof of Theorem~\ref{PositiveETS-Boundary-APS} then follows the same strategy as in the previous section, with a few additional points requiring clarification. We define $\mathbf{H}^\pm$ as the $||\,.\,||-$completion of $\widetilde{C}^{1,\pm}_c$. Then, for any $\chi\in L^2$, there exists a unique $\xi\in\mathbf{H}^+$ such that the spinor field $\xi_0:=\widetilde{D}^+\xi-\chi\in L^2$ is a weak solution of the system
	$$
	\left\lbrace
	\begin{array}{ll}
		\widetilde{D}^-\xi_0 =  0 \quad& {\rm on } \,\,M,\nonumber\\
		\mathcal{P}_{+}^-\xi_{0|\partial M}=0 \quad& {\rm along }\,\,\partial M\nonumber.
	\end{array}
	\right.
	$$
	Ellipticity of the Dirac-type operator $\widetilde{D}^-$ and the boundary condition $\mathcal{P}_{+}^-$ implies that $\xi_0\in\mathbf{H}^-\cap L^2$ is a strong solution of this boundary value problem. It follows from~(\ref{UnicityStep}) that $(\widetilde{D}^+)^k\xi_0\in L^2$ for all $k\in\mathbb{N}$, hence
	\begin{align*}
		\int_M|\widetilde{D}^+\xi_0|^2d\mu=\int_M\<\widetilde{D}^-\widetilde{D}^+\xi_0,\xi_0\>d\mu-in\oint_{\partial M}\<c_g(\nu)\xi_0,\xi_0\>d\sigma=-in\oint_{\partial M}\<c_g(\nu)\xi_0,\xi_0\>d\sigma
	\end{align*}
	since the operators $\widetilde{D}^+$ and $\widetilde{D}^+$ commute. Applying Lemma~\ref{LemmaAPS} to this identity yields  $\xi_0=0$. As a consequence, the operator $\widetilde{D}^+:\mathbf{H}^+\to L^2$ is an isomorphism, and the arguments of Section~\ref{BoundaryLessCase} apply verbatim, which completes the proof. The final ingredient is the following lemma.
	\begin{lemma}\label{LemmaAPS}
		For all $\xi_0\in\mathbf{H}^-$, we have:
		\begin{align*}
			-i\oint_{\partial M}\<c_g(\nu)\xi_0,\xi_0\>d\sigma\leq 0. 
		\end{align*}
	\end{lemma}
	
	{\it Proof.} We first claim that
	\begin{align*}
		-i\oint_{\partial M}\<c_g(\nu)\psi_k^-,\psi_l^-\>d\sigma=\frac{n-1}{\mu_k+\mu_l}\delta_{kl}
	\end{align*}
	for any two elements $\psi_k^-,\psi_l^-$ of an $L^2$-orthonormal basis of $\D^-$-eigenspinors. First, note that
	\begin{align*}
		-i \mu_k\oint_{\partial M}\<c_g(\nu)\psi_k^-,\psi_l^-\>d\sigma = i\oint_{\partial M}\<\D^+\left(c_g(\nu)\psi_k^-\right),\psi_l\>d\sigma
	\end{align*}
	by (\ref{RelationDpm}). Using the symmetry of $\D^+$, we obtain
	\begin{align*}
		-i \mu_k\oint_{\partial M}\<c_g(\nu)\psi_k^-,\psi_l^-\>d\sigma= i\oint_{\partial M}\<c_g(\nu)\psi_k^-,\D^+\psi_l^-\>d\sigma.
	\end{align*}
	On the other hand, by definition,
	\begin{align*}
		\D^+\psi= \D^-\psi+i(n-1)c_g(\nu)\psi
	\end{align*}
	for all $\psi\in\Gamma(\ETSB_g)$ and then 
	\begin{align*}
		-i \mu_k\oint_{\partial M}\<c_g(\nu)\psi_k^-,\psi_l^-\>d\sigma=i\mu_l\oint_{\partial M}\<c_g(\nu)\psi_k^+,\psi_l^-\>d\sigma+(n-1)\oint_{\partial M}\<\psi_k,\psi_l\>d\sigma.
	\end{align*}
	This proves the claim. Let now $\xi_0\in\mathbf{H}^-$. Writing its boundary trace as $\xi_0=\sum_{k<0}A_k\psi_k^-$, we get
	\begin{align*}
		-i\oint_{\partial M}\<c_g(\nu)\xi_0,\xi_0\>d\sigma = -i\sum_{k<0,\,l<0}A_k\overline{A_l}\oint_{\partial M}\<c_g(\nu)\psi_k^+,\psi_l^-\>d\sigma=\frac{n-1}{2}\sum_{k<0}\mu_k^{-1}|A_k|^2
	\end{align*}
	which is non-positive since $\mu_k<0$ for all $k<0$. 
	\qed
	
	This completes the proof of Theorem~\ref{PositiveETS-Boundary-APS}, by the same argument as in Section~\ref{BoundaryLessCase}.
	
	\subsection{The time-symmetric case: proof of Theorems~\ref{PMTC-Hyp} and ~\ref{PMTCB-Hyp-TS}}
	
	In the time-symmetric setting, the charged energy--momentum simplifies to
	\begin{align*}
		\widetilde{\Xi}(V,f):=\Xi(V,0,f)\quad\text{for}\quad(V,f)\in\mathcal{N}_b\oplus\mathbb{R}.
	\end{align*}
    One computes that
	\begin{align*}
		V_u=|\zeta_u|^2_b = 2\Big(|u|_\delta^2 V_{(0)}-\sum_{j=1}^n\<u,ic_\delta(e_j)u\>_\delta V_{(j)}\Big).
	\end{align*}
	Here, for $x\in\mathbb{B}^n$, the functions
	\begin{align*}
		V_{(0)}=\frac{1+|x|^2_\delta}{1-|x|^2_\delta}\quad\text{and}\quad V_{(j)}=\frac{2x^j}{1-|x|^2_\delta},\quad j=1,\cdots,n
	\end{align*}
	form a basis of the space of hyperbolic static potentials $\mathcal{N}_b$ in the unit ball model of $\mathbb{H}^n$. Similarly, one finds
	\begin{align*}
		f_u=\<\gamma_b\zeta_u,\zeta_u\>_b=2 \<\gamma_\delta u,u\>_\delta.
	\end{align*}
	
	It then follows from Theorem~\ref{PET-Spherical} (or from Theorems~\ref{PositiveETS-Boundary} and~\ref{PositiveETS-Boundary-APS} in the presence of an inner boundary),
	together with the linearity of $\widetilde{\Xi}$ and the expressions \eqref{MassCoeffSpherical}, that
	\begin{align*}
		0\leq 	\widetilde{\Xi}(V_u,f_u)= 2\<u,\Big(m_{(0)}u-i\sum_{j=1}^n m_{(j)}c_\delta(e_j)u+Q\gamma_\delta u\Big)\>_\delta
	\end{align*}
	for all $u\in\Sigma_n$. It is well known that this condition is equivalent to $m\geq |Q|$, which proves Theorem~\ref{PMTC-Hyp} and Theorem~\ref{PMTCB-Hyp-TS}.


	\bibliographystyle{alpha}     
	\bibliography{BiblioHabilitation}


	\end{document}